\documentclass{qjmam}

\usepackage{amsmath}
\usepackage{amssymb}
\usepackage{graphicx}
\usepackage{amsfonts}

\startpage{1}
\yr{2019}
\vol{XX}
\issue{X}

\makeatletter

    \makeatletter
    \renewenvironment{thebibliography}[1]
          {
           \@mkboth{\MakeUppercase\refname}{\MakeUppercase\refname}%
           \list{\@biblabel{\@arabic\c@enumiv}}%
                {\settowidth\labelwidth{\@biblabel{#1}}%
                 \leftmargin\labelwidth
                 \advance\leftmargin\labelsep
                 \@openbib@code
                 \usecounter{enumiv}%
                 \let\p@enumiv\@empty
                 \renewcommand\theenumiv{\@arabic\c@enumiv}}%
           \sloppy
           \clubpenalty4000
           \@clubpenalty \clubpenalty
           \widowpenalty4000%
           \sfcode`\.\@m}
          {\def\@noitemerr
            {\@latex@warning{Empty `thebibliography' environment}}%
           \endlist}
    \makeatother

\DeclareMathAlphabet\mathbit
    \encodingdefault\rmdefault\bfdefault\itdefault
\DeclareOldFontCommand{\bi}{\normalfont\bfseries\itshape}{\mathbit}

\makeatother

\begin{document}

\numberwithin{equation}{section}

\title[{Displacements~representations~for~the~problems}] {Displacements representations for the problems with spherical
and circular material surfaces with surface tension}

\author[S.~G.~MOGILEVSKAYA \etal] {\sc Sofia G. Mogilevskaya \and Volodymyr I. Kushch \and Anna Y. Zemlyanova}

\address{Department of Civil, Environmental,\\
and Geo- Engineering, University of Minnesota,\\ 
500 Pillsbury Drive S.E.,\\ 
Minneapolis, MN, 55455, USA }

\extraaddress{Institute for Superhard Materials\\ 
of the National Academy of Sciences of Ukraine,\\
04074 Kiev, Ukraine}

\extraaddress{Department of of Mathematics,\\
Kansas State University,\\
138 Cardwell Hall,\\
Manhattan, Kansas, 66506, USA}

\received{\recd XX.XX.XX. \revd XX.XX.XX}

\maketitle

\begin{abstract}
The displacements representations of the type used by Christensen
and Lo (1979) are modified to allow for analytical treatment of problems
involving spherical and circular material surfaces that possess constant
surface tension. The modified representations are used to derive closed-form
expressions for the local elastic fields and effective moduli of a
macroscopically isotropic composite materials containing spherical
and circular inhomogeneities with the interfaces described by the
complete Gurtin-Murdoch and Steigmann-Ogden models.
\end{abstract}


\section{Introduction}

In their paper, Christensen and Lo (1979) presented new micromechanical
scheme, the generalized self-consistent scheme, for the evaluation
of the effective shear moduli of macroscopically isotropic composites
containing cylindrical and spherical inhomogeneities. The scheme utilizes
the three phases model that involves a spherical or circular inhomogeneity,
a spherical or circular annulus of matrix material, and an infinite
outer region of equivalent homogeneous material subjected to uniform
far-field load.

For the case of arbitrary load, the rigorous solution of such model
problem would involve the use of infinite series (inner and outer
spherical harmonics in three dimensions or Laurent and Taylor series
in two dimensions) and could be quite complex. Instead, Christensen
and Lo (1979) suggested to use the closed-form representations for
the displacements in spherical and cylindrical coordinate systems
for the special far-field load - simple shear. The representations
were not new; they have been used earlier for various problems involving
spherical or circular boundaries, see Love (1927) for the three-dimensional
case and Savin (1961) for the two-dimensional one. However, they became
much more popular after publication of Christensen and Lo (1979) paper
and have been extensively used for the problems involving inhomogeneities
with interphases, imperfect interfaces, and material surfaces, see
e.g. Benveniste et al. (1989), Herve and Zaoui (1993), Duan et al.
(2006), Xu et al. (2016), Zemlyanova and Mogilevskaya (2018a), among
many others. 

We note that another type of displacement representations for the
axisymmetric problems with spherical surfaces was suggested in e.g.
Goodier (1933), see also Lurie (1964), and used for the problems with
imperfect interfaces and material surfaces in e.g. Duan et al. (2005,
2007), He and Li (2006), Lim et al. (2006), Mi and Kouris (2006).
While the solution for the case of simple shear far-field load can
be obtained from the superposition of the appropriately chosen axisymmetric
solutions, see e.g. Hashin (1991), the related procedure requires
an additional step and the derivations associated with that step. 

One limitation of Christensen and Lo (1979) type representations is
that they are not valid for the problems with surface and interface
effects as they do not include surface tension. Another limitation
is that they involve unknown coefficients, which in many cases have
to be found numerically from the system of linear algebraic equations,
thus not allowing for accurate identification of all governing problem
parameters. 

The concept of surface tension is important in modeling of various
nano-sized phenomena where the influence of the surface becomes more
significant due to the high surface-area-to-volume ratio, see e.g.
Cahn and Larché (1982), Miller and Shenoy (2000), Sharma and Ganti
(2002, 2004), Sharma et al. (2003), Dingreville et al. (2005), Lim
et al. (2006), He and Li (2006), Huang and Wang (2006), Mi and Kouris
(2006), Mogilevskaya et al. (2008, 2010), Ru (2010), Chhapadia et
al (2011), Kushch et al. (2011, 2013), Chatzigeorgiou et al. (2017),
Javili et al. (2017), and many others. 

The models proposed by Gurtin-Murdoch (1975, 1978) and Steigmann-Ogden
(1997, 1999) are the most studied continuum models of material surfaces
with surface tension. In these models, the interface between the material
constituents is treated as either a membrane (Gurtin-Murdoch model)
or a shell (Steigmann-Ogden model) of vanishing thickness possessing
surface tension as well as corresponding surface elastic properties.

Only a few closed-form solutions of the inhomogeneity problems involving
the \textit{complete} Gurtin-Murdoch model have been reported so far,
perhaps due notorious difficulties in handling effects introduced
by surface tension. Quite naturally, the benchmark problems, for which
analytical solutions could be constructed, involve regular shapes
of material surfaces. These solutions have been reported for the case
of circular inhomogeneity in e.g. Mogilevskaya et all. (2008), Jammes
et al. (2009), for elliptical inhomogeneity in Kushch et al. (2014),
for spherical inhomogeneity in Lim et al. (2006), He and Li (2006),
Kushch et al. (2011, 2013), and spheroidal inhomogeneity in Kushch
(2018), Kushch et al. (2018). The methods used in these publications
were quite complex, which was justified for the problems involving
multiple inhomogeneities. However, for the important case of a single
inhomogeneity, often used in various single-inhomogeneity-based homogenization
schemes, it is desirable to devise simpler solution by modifying the
representations used by Christensen and Lo (1979). 

To the best knowledge of the authors, no publication reports the solution
of the problem of a spherical inhomogeneity with the interface described
by the \textit{complete} Steigmann-Ogden model (with a full set of
interface parameters) and non-hydrostatic loading conditions. The
first and the only publication that considered the Steigmann-Ogden
model for spherical inhomogeneity was by Zemlyanova and Mogilevskaya
(2018a) who presented the solutions for two \textit{particular} cases:
(i) hydrostatic load (when the problem becomes one-dimensional) and
(ii) deviatoric load \& zero surface tension (when the problem can
be solved by using classical Christensen and Lo type displacement
representations). As stated above, the latter representations are
not valid for the cases involving surface tension and that was the
reason why the authors of that work have not been able to solve the
problem in more general setting. The two-dimensional solutions for
the problem of circular inhomogeneities with the complete Steigmann-Ogden
interface model are reported in Zemlyanova and Mogilevskaya (2018b)
and Han et al. (2018). However, they too were obtained with tedious
algebra that, for the case of a single inhomogeneity, could be avoided,
if simpler representations of the type used by Christensen and Lo
(1979) could be modified to include surface tension. 

Thus, the goal of the present paper is three-fold. First, we present
\textit{new} Christensen and Lo (1979) type \textit{representations}
to allow for simple analytical solutions of the problems involving
spherical and circular material surfaces possessing constant surface
tension. Second, using the obtained representations, we derive \textit{new
analytical solution }for the problem of a spherical inhomogeneity
with the interface described by the \textit{complete} Steigmann-Ogden
model. Finally, we provide \textit{closed-form analytical expressions}
for the coefficients involved in those representations for the cases
of the complete Gurtin-Murdoch and Steigmann-Ogden models. We emphasize
again that the focus of the present paper is not on the modification
of micromechanical scheme of Christensen and Lo (1979) but rather
on the representations they used to obtain it. The new representations
devised here could be used in variety of micromechanical schemes one
of which, Maxwell's (1873) homogenization scheme, is used in the present
paper.

The paper is structured as follows. In Section 2, we review the Christensen
and Lo representations for the case of simple shear far-field load
and, for completeness, provide similar representation for the case
of hydrostatic load. In Section 3 (with the details provided in Appendix
A), we use analytical solutions for the problems of a circular inhomogeneity
with the complete Gurtin-Murdoch (Mogilevskaya et al., 2008) and Steigmann-Ogden
(Zemlyanova and Mogilevskaya, 2018b) interfaces to construct the Christensen
and Lo (1979) type representations for the hydrostatic and simple
shear far-field loads. In Section 4 (with the details provided in
Appendix B), we use analytical solution of Kushch et al. (2011) to
presents analogous representations for the problem of a spherical
inhomogeneity with the complete Gurtin-Murdoch interface. In Section
5, we use the representations of Section 4 to obtain new analytical
solution for the problem of a spherical inhomogeneity with the interface
described by the complete Steigmann-Ogden model. In Section 6, we
use this solution in combination with Maxwell's (1873) methodology
to derive the single-inhomogeneity based estimate of the effective
shear modulus of macroscopically isotropic material containing spherical
inhomogeneities with interfaces described by the complete Steigmann-Ogden
model. In Section 7, we present a summary of our results and conclusions.

\smallskip{}

\section{Christensen and Lo (1979) representations for
circular and spherical inhomogeneities}

Consider the problem of a circular (Fig.1a) or a spherical (Fig.1b)
inhomogeneity of radius $R$ embedded into an infinite matrix and
subjected to uniform far field load $\sigma_{kj}^{\infty}$, $k,j=1,\ldots,d$,
where $d=2$ in two dimensions and $d=3$ in three dimensions. The
center of the inhomogeneity is located at the origin of the Cartesian
coordinate system. Assume also that the bulk material of the matrix
(inhomogeneity) is linearly elastic and isotropic; the corresponding
elastic moduli for the two phases are shear modulus $\mu$($\mu_{I}$)
and Poisson's ratio $\nu$($\nu_{I}$). The interface conditions between
the matrix and inhomogeneity are not specified at this time.

\begin{figure}[ht]
\centering\includegraphics[width=0.8\textwidth]{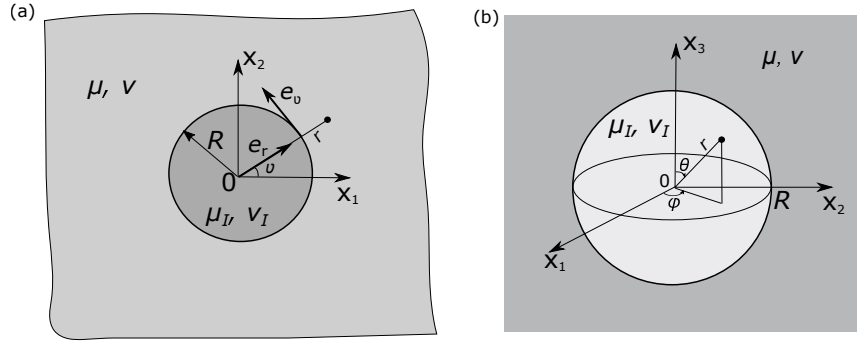} \caption{a) Circular and b) spherical inhomogeneity in an infinite matrix}
\label{fig_1.1} 
\end{figure}

Christensen and Lo (1979) used the following representations for the
displacements in the polar or spherical coordinate systems of Fig.1
and the case of simple shear far-field load (the only non-vanishing
far-field stresses are $\sigma_{11}^{\infty}=-\sigma_{22}^{\infty}=\sigma_{d}^{\infty}$):
\begin{itemize}
\item circular inhomogeneity
\end{itemize}
\textit{inside the inhomogeneity}

\begin{align}
u_{r}(z) & =\frac{R}{4\mu_{I}}\left[d_{1}\frac{r}{R}+\left(\kappa_{I}-3\right)a_{1}\frac{r^{3}}{R^{3}}\right]\cos2\vartheta\label{eq:2.1}\\
u_{\vartheta}(z) & =\frac{R}{4\mu_{I}}\left[-d_{1}\frac{r}{R}+\left(\kappa_{I}+3\right)a_{1}\frac{r^{3}}{R^{3}}\right]\sin2\vartheta\nonumber 
\end{align}

\textit{inside the matrix }

\begin{align}
u_{r}(z) & =\frac{R}{4\mu}\left[2\sigma_{d}^{\infty}\frac{r}{R}+\left(\kappa+1\right)a_{3}\frac{R}{r}+c_{3}\frac{R^{3}}{r^{3}}\right]\cos2\vartheta\label{eq:2.2}\\
u_{\vartheta}(z) & =\frac{R}{4\mu}\left[-2\sigma_{d}^{\infty}\frac{r}{R}-\left(\kappa-1\right)a_{3}\frac{R}{r}+c_{3}\frac{R^{3}}{r^{3}}\right]\sin2\vartheta\nonumber 
\end{align}
in which $\kappa=3-4\nu,\:\kappa_{I}=3-4\nu_{I}$, and the unknown
coefficients $d_{1},\,a_{1},\,a_{3},\,c_{3}$ are found from the interface
and far-field conditions. It should be noted that the expressions
of Eq. (4.1) in Christensen and Lo (1979) were dimensionally inconsistent,
so in Eq. (\ref{eq:2.2}) we added missing multipliers to eliminate
this inconsistency.
\begin{itemize}
\item spherical inhomogeneity 
\end{itemize}
\begin{align}
u_{r} & =U_{r}\left(r\right)\sin^{2}\theta\cos2\varphi\label{eq:2.3}\\
u_{\theta} & =U_{\theta}\left(r\right)\sin\theta\cos\theta\cos2\varphi\nonumber \\
u_{\varphi} & =-U_{\theta}\left(r\right)\sin\theta\sin2\varphi\nonumber 
\end{align}
in which the functions $U_{r}\left(r\right)$, $U_{\theta}(r)$ of
Eq. (\ref{eq:2.3}) are different for the inhomogeneity and the matrix
and taken as 

\textit{inside the inhomogeneity }

\begin{align}
U_{r}^{inh}\left(r\right) & =A_{1}r-\frac{6\nu_{I}}{1-2\nu_{I}}A_{2}r^{3}\label{eq:2.4}\\
U_{\theta}^{inh}\left(r\right) & =A_{1}r-\frac{7-4\nu_{I}}{1-2\nu_{I}}A_{2}r^{3}\nonumber 
\end{align}

\textit{inside the matrix }

\begin{align}
U_{r}^{mat}\left(r\right) & =D_{1}r+\frac{3D_{3}}{r^{4}}+\frac{5-4\nu}{1-2\nu}\frac{D_{4}}{r^{2}}\label{eq:2.5}\\
U_{\theta}^{mat}\left(r\right) & =D_{1}r-\frac{2D_{3}}{r^{4}}+2\frac{D_{4}}{r^{2}}\nonumber 
\end{align}
and involve five unknown coefficients $A_{1},\:A_{2},\:D_{1},\:D_{3},\:D_{4}$
that have to be found from the interface and far-field conditions. 

Note that the representations of Eqs. (\ref{eq:2.1})-(\ref{eq:2.5})
are not valid for the problems with surface tension. We add also that
Christensen and Lo (1979) considered a three-phase model and their
set of representations included expressions for the displacements
inside the third phase (interphase), which is not present in the problems
under study that deal with material surfaces of vanishing thicknesses.

For the completeness, we also list here the representations for the
case of hydrostatic load (the only non-vanishing far-field stresses
are $\sigma_{kk}^{\infty}=\sigma_{h}^{\infty},\:k=1,\ldots,d$ ).
Such representations are often used for the evaluation of the effective
bulk modulus. The problem is one-dimensional and the expressions for
the only non-vanishing radial component of the displacements are
\begin{itemize}
\item inside the inhomogeneity 
\end{itemize}
\begin{align}
u_{r}(z) & =F_{1}r\label{eq:2.6}
\end{align}

\begin{itemize}
\item inside the matrix
\end{itemize}
\begin{align}
u_{r}(z) & =F_{2}r+F_{3}/r^{d-1}\label{eq:2.7}
\end{align}
where again $d$ is the dimension of the problem and the three unknown
coefficients $F_{1},\:F_{2},\:F_{3}$ have to be found from the interface
and far-field conditions.

\section{Representations for a circular inhomogeneity
with the Gurtin-Murdoch and Steigmann-Ogden interfaces}

Consider circular inhomogeneity shown on Fig.1a. In both models, the
displacements are continuous across the interface but the tractions
undergo jumps. The tractions jump conditions for the more general
Steigmann-Ogden model (characterized by the elastic parameters $\mu_{0}$,
$\lambda_{0}$, surface tension $\sigma_{0}$, and bending parameters
$\zeta_{0}$ and $\chi_{0}$) can be written as, see Zemlyanova and
Mogilevskaya (2018b)

\begin{align}
\sigma_{rr}^{inh}-\sigma_{rr}^{mat} & =-\frac{\sigma_{0}}{R}+\frac{\sigma_{0}}{R^{2}}\left(u_{r,\vartheta\vartheta}-u_{\vartheta,\vartheta}\right)\label{eq:3.1}\\
- & \left(\lambda_{0}+2\mu_{0}\right)\frac{1}{R^{2}}\left(u_{\vartheta,\vartheta}+u_{r}\right)\nonumber \\
- & \left(2\chi_{0}+\zeta_{0}\right)\frac{1}{R^{4}}\left(u_{r,\vartheta\vartheta\vartheta\vartheta}-u_{\vartheta,\vartheta\vartheta\vartheta}\right)\nonumber 
\end{align}

\begin{align}
\sigma_{r\vartheta}^{inh}-\sigma_{r\vartheta}^{mat} & =\frac{\sigma_{0}}{R^{2}}\left(u_{r,\vartheta}-u_{\vartheta}\right)\label{eq:3.2}\\
+ & \left(\lambda_{0}+2\mu_{0}\right)\frac{1}{R^{2}}\left(u_{\vartheta,\vartheta\vartheta}+u_{r,\vartheta}\right)\nonumber \\
- & \left(2\chi_{0}+\zeta_{0}\right)\frac{1}{R^{4}}\left(u_{r,\vartheta\vartheta\vartheta}-u_{\vartheta,\vartheta\vartheta}\right)\nonumber 
\end{align}
in which $u_{r}$ and $u_{\vartheta}$ are local components of the
surface displacements in the local coordinate system ($r,\vartheta)$
shown in Fig 1a, $\sigma_{rr},\:\sigma_{r\vartheta}$ are the corresponding
components of tractions in that system, and the subscript ``,''
indicates differentiation, e.g. $u_{\vartheta,\vartheta}=\partial u_{\vartheta}/\partial\vartheta$. 

The tractions jump conditions for the Gurtin-Murdoch model could be
recovered from Eqs. (\ref{eq:3.1}), (\ref{eq:3.2}) by assuming that
the bending parameters vanish, i.e. $\zeta_{0}=0$ and $\chi_{0}=0$.

If the entire system is subjected to uniform far-field load, both
problems can be solved analytically using complex variables formalism,
see Mogilevskaya et al. (2008), Zemlyanova and Mogilevskaya (2018b).
Using these solutions (briefly reviewed in Appendix A) and some algebra,
the polar coordinates of the displacements everywhere in the composite
system can be expressed in the following closed forms: 

\textit{inside the inhomogeneity}

\begin{align}
u_{r} & =-\frac{\sigma_{0}}{2\left(K_{I}^{\left(2\right)}+\mu+2\eta\right)}\frac{r}{R}+\frac{R}{4\mu_{I}}\left[d_{1}\frac{r}{R}+\left(\kappa_{I}-3\right)a_{1}\frac{r^{3}}{R^{3}}\right]\cos2\vartheta\label{eq:3.21}\\
u_{\vartheta} & =\frac{R}{4\mu_{I}}\left[-d_{1}\frac{r}{R}+\left(\kappa_{I}+3\right)a_{1}\frac{r^{3}}{R^{3}}\right]\sin2\vartheta\nonumber 
\end{align}

\textit{inside the matrix}

\begin{align}
u_{r} & =-\frac{\sigma_{0}}{2\left(K_{I}^{\left(2\right)}+\mu+2\eta\right)}\frac{R}{r}+\frac{R}{4\mu}\left[2\sigma_{d}^{\infty}\frac{r}{R}+\left(\kappa+1\right)a_{3}\frac{R}{r}+c_{3}\frac{R^{3}}{r^{3}}\right]\cos2\vartheta\label{eq:3.22}\\
u_{\vartheta} & =\frac{R}{4\mu}\left[-2\sigma_{d}^{\infty}\frac{r}{R}-\left(\kappa-1\right)a_{3}\frac{R}{r}+c_{3}\frac{R^{3}}{r^{3}}\right]\sin2\vartheta\nonumber 
\end{align}
in which $K_{I}^{\left(2\right)}=2\mu_{I}/\left(\kappa_{I}-1\right)$
is the two-dimensional bulk modulus of the inhomogeneity, 
\begin{align}
d_{1} & =\frac{4\mu_{I}}{\kappa_{I}}\left(3\frac{A_{3}}{R}+\kappa_{I}\frac{A_{-1}}{R}\right)\label{eq:3.25}\\
a_{1} & =\frac{4\mu_{I}}{\kappa_{I}}\frac{A_{3}}{R}\nonumber \\
a_{3} & =\frac{4}{\kappa+1}\left(3\eta^{\left(2\right)}\frac{A_{3}}{R}-\omega^{\left(1\right)}\frac{A_{-1}}{R}\right)\nonumber \\
c_{3} & =\frac{4}{\kappa+1}\left[\left(\omega^{\left(1\right)}+\kappa\eta^{\left(2\right)}\right)\frac{A_{-1}}{R}-\omega^{\left(2\right)}\frac{A_{3}}{R}\right]\nonumber 
\end{align}
and the coefficients $A_{-1},\:A_{3}$ are given, for the case of
simple shear far-field load, by the last two expressions of Eqs. (\ref{eq:3.20})
of Appendix A. The meanings of remaining parameters involved in Eqs.
(\ref{eq:3.21})-(\ref{eq:3.25}) are also explained in that appendix.

It follows from the representations of Eqs. (\ref{eq:3.21}), (\ref{eq:3.22})
that the only difference between them and those of Christensen and
Lo (1979) are in radial components of the displacements, which now
has the following forms:

\begin{equation}
u_{r}=u_{r}^{C-L}+A\sigma_0\left(\frac{r}{R}\right)^{p}\label{eq:3.24}
\end{equation}
where $A=-\frac{1}{2\left(K_{I}^{\left(2\right)}+\mu+2\eta\right)}$, $p=1$ inside the inhomogeneity, $p=-1$ inside the matrix and
$u_{r}^{C-L}$ are given by Eqs. (\ref{eq:2.1}), (\ref{eq:2.2}). 

From Eq. (\ref{eq:A.6}) of Appendix A, it follows that, for the case
of hydrostatic far-field load $\sigma_{11}^{\infty}=\sigma_{22}^{\infty}=\sigma_{h}^{\infty},\:\sigma_{12}^{\infty}=0$,
the coefficients involved in Eqs. (\ref{eq:3.10}), (\ref{eq:3.11})
of that appendix are
\begin{equation}
ReA_{1}=\frac{R}{2\left(K_{I}^{\left(2\right)}+\mu+2\eta\right)}\left(\frac{\kappa+1}{2}\sigma_{h}^{\infty}-\frac{\sigma_{0}}{R}\right),\:A_{-1}=A_{3}=0\label{eq:3.19}
\end{equation}

Substitution of the coefficients of Eqs. (\ref{eq:3.19}) into Eqs.
(\ref{eq:3.10}), (\ref{eq:3.11}) leads to the representations for
the displacements given by Eqs. (\ref{eq:2.6}), (\ref{eq:2.7}) with
$d=2$ and the following coefficients $F_{1},\:F_{2},\:F_{3}$:

\begin{align}
F_{1} & =\frac{1}{2\left(K_{I}^{\left(2\right)}+\mu+2\eta\right)}\left(\frac{\kappa+1}{2}\sigma_{h}^{\infty}-\frac{\sigma_{0}}{R}\right)\label{eq:3.23}\\
F_{2} & =\frac{\sigma_{h}^{\infty}}{2K^{\left(2\right)}}\nonumber \\
F_{3} & =-\frac{R^{2}}{2\left(K_{I}^{\left(2\right)}+\mu+2\eta\right)}\left[\sigma_{h}^{\infty}\left(\frac{K_{I}^{\left(2\right)}}{K^{\left(2\right)}}-1+\frac{2\eta}{K^{\left(2\right)}}\right)+\frac{\sigma_{0}}{R}\right]\nonumber 
\end{align}

where $K^{\left(2\right)}=2\mu/\left(\kappa-1\right)$ is two-dimensional
bulk modulus of the matrix.

\section{Representations for a spherical inhomogeneity
with the Gurtin-Murdoch material surface}

Consider spherical inhomogeneity shown on Fig.1b. In\textbf{ }the
Gurtin-Murdoch model, the displacements are continuous across the
interface but the tractions undergo jumps. The tractions jump conditions
for the model (in which the interface is characterized by the elastic
parameters $\mu_{0}$, $\lambda_{0}$ and the surface tension $\sigma_{0}$)
can be written as, see Kushch et al. (2013), Zemlyanova and Mogilevskaya
(2018a), (the latter paper had a misprint, missing term $T_{\theta r}$
in Eq. (\ref{eq:4.2}))

\begin{equation}
\sigma_{rr}^{inh}-\sigma_{rr}^{mat}=\frac{1}{r\sin\theta}\left[T_{\varphi r,\varphi}+T_{\theta r}\cos\theta+(T_{\theta r,\theta}-T_{\varphi\varphi}-T_{\theta\theta})\sin\theta\right]\label{eq:4.1}
\end{equation}
\begin{equation}
\sigma_{r\theta}^{inh}-\sigma_{r\theta}^{mat}=\frac{1}{r\sin\theta}\left[T_{\varphi\theta,\varphi}+\left(T_{\theta\theta,\theta}+T_{\theta r}\right)\sin\theta+(T_{\theta\theta}-T_{\varphi\varphi})\cos\theta\right]\label{eq:4.2}
\end{equation}
\begin{equation}
\sigma_{r\varphi}^{inh}-\sigma_{r\varphi}^{mat}=\frac{1}{r\sin\theta}\left[T_{\varphi\varphi,\varphi}+(T_{\theta\varphi,\theta}+T_{\varphi r})\sin\theta+(T_{\varphi\theta}+T_{\theta\varphi})\cos\theta\right]\label{eq:4.3}
\end{equation}
in which 

\begin{equation}
T_{\theta\theta}=\sigma_{0}+\frac{\lambda_{0}+\sigma_{0}}{r\sin\theta}(u_{\varphi,\varphi}+u_{\theta}\cos\theta+u_{r}\sin\theta)+\frac{\lambda_{0}+2\mu_{0}}{r}(u_{\theta,\theta}+u_{r})\label{eq:4.4}
\end{equation}

\begin{equation}
T_{\varphi\varphi}=\sigma_{0}+\frac{\lambda_{0}+2\mu_{0}}{r\sin\theta}(u_{\varphi,\varphi}+u_{\theta}\cos\theta+u_{r}\sin\theta)+\frac{\lambda_{0}+\sigma_{0}}{r}(u_{\theta,\theta}+u_{r})\label{eq:4.5}
\end{equation}
\begin{equation}
T_{\varphi\theta}=\frac{1}{r\sin\theta}\left[\mu_{0}(u_{\theta,\varphi}-u_{\varphi}\cos\theta)+(\mu_{0}-\sigma_{0})u_{\varphi,\theta}\sin\theta\right]\label{eq:4.6}
\end{equation}
\begin{equation}
T_{\theta\varphi}=\frac{1}{r\sin\theta}\left[(\mu_{0}-\sigma_{0})(u_{\theta,\varphi}-u_{\varphi}\cos\theta)+\mu_{0}u_{\varphi,\theta}\sin\theta\right]\label{eq:4.7}
\end{equation}
\begin{equation}
T_{\theta r}=\frac{\sigma_{0}}{r}(u_{r,\theta}-u_{\theta})\label{eq:4.8}
\end{equation}
\begin{equation}
T_{\varphi r}=\frac{\sigma_{0}}{r\sin\theta}(u_{r,\varphi}-u_{\varphi}\sin\theta)\label{eq:4.9}
\end{equation}

If the entire system is subjected to uniform far-field load, the problem
can be solved analytically using the technique of vector spherical
harmonics, see Kushch (2013), Kushch et al. (2011). In that technique,
the displacement vector in spherical coordinates is sought as a linear
combination of vector partial solutions of the Lamé equation of isotropic
elasticity and the coefficients involved in the combination are found
from the boundary conditions.

For the case of simple shear far-field load, the only non-zero components
are $\sigma_{11}^{\infty}=-\sigma_{22}^{\infty}=\sigma_{d}^{\infty}$,
the displacement vector fields $\mathbf{u}\left(u_{r},u_{\theta},u_{\varphi}\right)$
can be represented inside the inhomogeneity and matrix as

\textit{inside the inhomogeneity }($\left\vert \mathbf{x}\right\vert =r<R$)
\begin{equation}
\mathbf{u}(\mathbf{x})=\frac{3}{2(1-2\nu_{I})}A_{0}\mathbf{u}_{00}^{(3)}(\mathbf{x})+4A_{1}Re\mathbf{u}_{22}^{(1)}(\mathbf{x})-\frac{84}{1-2\nu_{I}}A_{2}Re\mathbf{u}_{22}^{(3)}(\mathbf{x}),\label{eq:4.10}
\end{equation}

\textit{inside the matrix }($\left\vert \mathbf{x}\right\vert =r>R$)
\begin{equation}
\mathbf{u}(\mathbf{x})=2\frac{\sigma_{d}^{\infty}}{\mu}Re\mathbf{u}_{22}^{(1)}(\mathbf{x})+D_{0}\mathbf{U}_{00}^{(1)}(\mathbf{x})-\frac{D_{3}}{3}Re\mathbf{U}_{22}^{(1)}(\mathbf{x})+\frac{D_{4}}{1-2\nu}Re\mathbf{U}_{22}^{(3)}(\mathbf{x}),\label{eq.:4.11}
\end{equation}
with $A_{0}$, $A_{1}$, $A_{2}$ and $D_{0}$, $D_{3}$ $D_{4}$
being the unknown constants and $\mathbf{u}_{kk}^{(m)}(\mathbf{x}),\mathbf{U}_{kk}^{(m)}(\mathbf{x})$
being the vector partial solutions of the Lame equation in spherical
coordinates. These functions are given by Eqs. (\ref{eq:B.1}), (\ref{eq:B.2})
of Appendix B. 

After some algebra that involves the use of Eqs. (\ref{eq:B.1})-(\ref{eq:B.3})
of Appendix B, one arrives at the following component representations
of Eqs. (\ref{eq:4.10}), (\ref{eq.:4.11}):

\textit{inside the inhomogeneity} 
\begin{align}
u_{r} & =-rA_{0}+U_{r}^{inh}\left(r\right)\sin^{2}\theta\cos2\varphi,\label{eq:4.12}\\
u_{\theta} & =U_{\theta}^{inh}\left(r\right)\sin\theta\cos\theta\cos2\varphi,\label{eq:4.13}\\
u_{\varphi} & =-U_{\theta}^{inh}\left(r\right)\sin\theta\sin2\varphi,\label{eq:4.14}
\end{align}

\textit{inside the matrix }
\begin{align}
u_{r} & =-\frac{1}{r^{2}}D_{0}+U_{r}^{mat}\left(r\right)\sin^{2}\theta\cos2\varphi\label{eq:4.15}\\
u_{\theta} & =U_{\theta}^{mat}\left(r\right)\sin\theta\cos\theta\cos2\varphi,\label{eq:4.16}\\
u_{\varphi} & =-U_{\theta}^{mat}\left(r\right)\sin\theta\sin2\varphi,\label{eq:4.17}
\end{align}
in which $U_{r}^{inh}\left(r\right)$, $U_{r}^{mat}\left(r\right)$
are the functions defined in Eqs. (\ref{eq:2.4}), (\ref{eq:2.5})
with $D_{1}=\sigma_{d}^{\infty}/2\mu$. The remaining coefficients
involved in Eqs. (\ref{eq:4.12})-(\ref{eq:4.17}) can be found from
the system of Eqs. (\ref{eq:B.7})-(\ref{eq:B.9}) presented in Appendix
B.

It can be seen that the only difference between the representations
of Eqs. (\ref{eq:4.12})-(\ref{eq:4.17}) and those of Christensen
and Lo (1979) is in radial displacements, which now have the following
forms that involve additional terms due to surface tension:

\begin{equation}
u_{r}=u_{r}^{C-L}+A\sigma_0\left(\frac{r}{R}\right)^{p}\label{eq:4.18}
\end{equation}
where $A=-\frac{2}{4\mu+3K_{I}^{\left(3\right)}+2\eta_{0}}$, $K_{I}^{\left(3\right)}=\frac{2}{3}\mu_{I}\left(1+\nu_{I}\right)/\left(1-2\nu_{I}\right)$
is three-dimensional bulk modulus of the inhomogeneity, $p=1$ inside
the inhomogeneity, $p=-2$ inside the matrix, $u_{r}^{C-L}$ are given
by Eqs. (\ref{eq:2.3})- (\ref{eq:2.5}), and 
\begin{equation}
\eta_{0}=\left(2\mu_{0}+2\lambda_{0}+\sigma_{0}\right)/R\label{eq:4.23}
\end{equation}

For the case of hydrostatic far-field load ($\sigma_{kk}^{\infty}=\sigma_{h}^{\infty},\:\sigma_{kj}^{\infty}=0,\:k\neq j$),
the representations for the displacements are given by Eqs. (\ref{eq:2.6}),
(\ref{eq:2.7}) with $d=3$ and the following coefficients:

\begin{align}
F_{1} & =\frac{\left[1+4\mu/\left(3K^{\left(3\right)}\right)\right]\sigma_{h}^{\infty}-2\sigma_{0}/R}{4\mu+3K_{I}^{\left(3\right)}+2\eta_{0}}\label{eq:4.22}\\
F_{2} & =\frac{\sigma_{h}^{\infty}}{3K^{\left(3\right)}}\nonumber \\
F_{3} & =R^{3}\frac{\sigma_{h}^{\infty}\left[\left(1-K_{I}^{\left(3\right)}/K^{\left(3\right)}\right)-2\eta_{0}/\left(3K^{\left(3\right)}\right)\right]-2\sigma_{0}/R}{4\mu+3K_{I}^{\left(3\right)}+2\eta_{0}}\nonumber 
\end{align}
where $K^{\left(3\right)}=\frac{2}{3}\mu\left(1+\nu\right)/\left(1-2\nu\right)$
is three-dimensional bulk modulus of the matrix.\smallskip{}

\section{New solutions for a spherical inhomogeneity with
the Steigmann-Ogden material surface (simple shear far-field load)}

\noindent Consider again spherical inhomogeneity shown on Fig.1b but
now assume that its interface is described by the complete Steigmann-Ogden
model that includes surface tension. In this case, the jump of tractions
on the spherical boundary can be written as, see Zemlyanova and Mogilevskaya
(2018a): 
\begin{equation}
\sigma_{rr}^{inh}-\sigma_{rr}^{mat}=\frac{1}{r\sin\theta}\left[T_{\varphi r,\varphi}+T_{\theta r}\cos\theta+(T_{\theta r,\theta}-T_{\varphi\varphi}-T_{\theta\theta})\sin\theta\right]+\label{eq:5.1}
\end{equation}
\[
\frac{1}{r^{2}\sin^{2}\theta}\left[M_{\varphi\varphi,\varphi\varphi}+(M_{\theta\varphi,\theta\varphi}+M_{\varphi\theta,\varphi\theta})\sin\theta+M_{\theta\theta,\theta\theta}\sin^{2}\theta+(M_{\varphi\theta,\varphi}+M_{\theta\varphi,\varphi})\cos\theta+\right.
\]
\[
\left.(2M_{\theta\theta,\theta}-M_{\varphi\varphi,\theta})\cos\theta\sin\theta-(M_{\theta\theta}-M_{\varphi\varphi})\sin^{2}\theta\right]
\]
\begin{equation}
\sigma_{r\theta}^{inh}-\sigma_{r\theta}^{mat}=\frac{1}{r\sin\theta}\left[T_{\varphi\theta,\varphi}+(T_{\theta\theta,\theta}+T_{\theta r})\sin\theta+(T_{\theta\theta}-T_{\varphi\varphi})\cos\theta\right]+\label{eq:5.2}
\end{equation}
\[
\frac{1}{r^{2}\sin\theta}\left[M_{\varphi\theta,\varphi}+(M_{\theta\theta}-M_{\varphi\varphi})\cos\theta+M_{\theta\theta,\theta}\sin\theta\right]
\]
\begin{equation}
\sigma_{r\varphi}^{inh}-\sigma_{r\varphi}^{mat}=\frac{1}{r\sin\theta}\left[T_{\varphi\varphi,\varphi}+(T_{\theta\varphi,\theta}+T_{\varphi r})\sin\theta+(T_{\varphi\theta}+T_{\theta\varphi})\cos\theta\right]+\label{eq:5.3}
\end{equation}
\[
\frac{1}{r^{2}\sin\theta}\left[M_{\varphi\varphi,\varphi}+(M_{\varphi\theta}+M_{\theta\varphi})\cos\theta+M_{\theta\varphi,\theta}\sin\theta\right]
\]
in which the components of the surface stress tensor ${\bf T}$ are
given by Eqs. (\ref{eq:4.4})-(\ref{eq:4.9}), and the components
of the surface couple-stress tensor ${\bf M}$ are 
\begin{equation}
M_{\varphi\varphi}=(2\chi_{0}+\zeta_{0})\kappa_{\varphi\varphi}+\zeta_{0}\kappa_{\theta\theta}\label{eq:5.4}
\end{equation}
\begin{equation}
M_{\varphi\theta}=M_{\theta\varphi}=2\chi_{0}\kappa_{\varphi\theta}\label{eq:5.5}
\end{equation}
\begin{equation}
M_{\theta\theta}=\zeta_{0}\kappa_{\varphi\varphi}+(2\chi_{0}+\zeta_{0})\kappa_{\theta\theta}\label{eq:5.6}
\end{equation}
\begin{equation}
\kappa_{\varphi\varphi}=-\frac{1}{r^{2}\sin^{2}\theta}\left[u_{r,\varphi\varphi}-u_{\varphi,\varphi}\sin\theta+(u_{r,\theta}-u_{\theta})\cos\theta\sin\theta\right]\label{eq:5.7}
\end{equation}
\begin{equation}
\kappa_{\varphi\theta}=\kappa_{\theta\varphi}=\frac{1}{2r^{2}\sin^{2}\theta}\left[(u_{r,\varphi\theta}-u_{\theta,\varphi})\sin\theta-(u_{r,\varphi}-u_{\varphi}\sin\theta)\cos\theta+\frac{\partial}{\partial\theta}\left(\frac{u_{r,\varphi}}{\sin\theta}-u_{\varphi}\right)\sin^{2}\theta\right]\label{eq:5.8}
\end{equation}
\begin{equation}
\kappa_{\theta\theta}=\frac{1}{r^{2}}(u_{r,\theta\theta}-u_{\theta,\theta})\label{eq:5.9}
\end{equation}
Eqs. (\ref{eq:5.4})-(\ref{eq:5.9}) involve the components of the
tensor of changes of curvature $\boldsymbol{\kappa}$ and the bending
stiffness parameters $\chi_{0}$, $\zeta_{0}$.

Now, assuming simple shear load at infinity, we propose to use the
representations of Eqs. (\ref{eq:4.12})-(\ref{eq:4.17}) for obtaining
the new solution for the problem of a spherical inhomogeneity with
the complete Steigmann-Ogden interface. In Zemlyanova and Mogilevskaya
(2018a), the solution for the simple shear far-field load was only
presented for the case of zero surface tension.

Substitution of the representations of Eqs. (\ref{eq:4.12})-(\ref{eq:4.17})
into Eqs. (\ref{eq:5.1})-(\ref{eq:5.3}) leads to the following expressions
for the jumps of boundary tractions: 
\[
\sigma_{r\varphi}^{inh}-\sigma_{r\varphi}^{mat}=\left[\frac{2}{r^{2}}\left((3\lambda_{0}+5\mu_{0}+\sigma_{0})U_{\theta}-2(\mu_{0}+\lambda_{0}+\sigma_{0})U_{r}\right)-2\gamma\frac{R^{3}}{r^{4}}(2U_{r}-U_{\theta})\right]\sin\theta\sin2\varphi
\]
\[
\sigma_{r\theta}^{inh}-\sigma_{r\theta}^{mat}=
\]
\begin{equation}
\left[\frac{1}{r^{2}}\left(-(3\lambda_{0}+5\mu_{0}+\sigma_{0})U_{\theta}+2(\mu_{0}+\lambda_{0}+\sigma_{0})U_{r}\right)+\gamma\frac{R^{3}}{r^{4}}(2U_{r}-U_{\theta})\right]\sin2\theta\cos2\varphi\label{eq:5.10}
\end{equation}
\[
\sigma_{rr}^{inh}-\sigma_{rr}^{mat}=-\frac{2\sigma_{0}}{r}+
\]
\[
\left[\frac{1}{r^{2}}\left(6(\lambda_{0}+\mu_{0}+\sigma_{0})U_{\theta}-4(\mu_{0}+\lambda_{0}+2\sigma_{0})U_{r}\right)-6\gamma\frac{R^{3}}{r^{4}}(2U_{r}-U_{\theta})\right]\sin^{2}\theta\cos2\varphi
\]
where $\gamma=\left(3\chi_{0}+5\zeta_{0}\right)/R^{3}$.

The first and third equations of Eqs. (\ref{eq:5.10}) are identical
to the corresponding equations involved in Eqs. (B3) from Zemlyanova
and Mogilevskaya (2018a). The second equation of Eqs. (\ref{eq:5.10})
is corrected version of the second equation involved in Eqs. (B3)
from the same paper, which had a misprint that, however, have not
affected the results presented there.

The use of Eqs. (\ref{eq:5.10}) together with the conditions of continuity
of displacements across the boundary of the sphere leads to the linear
system of four equations to find the unknown coefficients Eqs.(\ref{eq:2.4})-(\ref{eq:2.5}).
This linear system and its solution are presented in the Appendix
C.

To illustrate the effects of the Steigmann-Ogden interface parameters,
we consider the example similar to that presented in Zemlyanova and
Mogilevskaya (2018a) that involves a cavity of radius $R=5\,nm$ and
assume that the normalized simple shear stress at infinity is 

\begin{equation}
\sigma_{d}/\mu=0.000028818\label{eq:5.11}
\end{equation}

We also adopt the following three values of the normalized surface
tension:

\begin{equation}
\sigma_{0}/\mu R=0,\:\sigma_{0}/\mu R=0.0067435,\:\sigma_{0}/\mu R=0.0097983\label{eq:5.12}
\end{equation}

The rest of the parameters for the example are chosen to be identical
to the following ones used in Zemlyanova and Mogilevskaya (2018a):

\begin{equation}
\nu=0.3,\:\mu_{0}/\mu R=0.030156,\:\lambda_{0}/\mu R=0.060312,\:\gamma/\mu=0.00028382\label{eq:5.13}
\end{equation}

On Fig. 2, we plotted the normalized hoop stress $\sigma_{\theta\theta}/\mu$
along the meridian line $\varphi=0$ at the cavity surface ($0\leq \theta \leq \pi/2$).
It can be seen that the surface tension has significant effect on
the normalized hoop stress variation. The normalized stress becomes
compressive when $\sigma_{0}\neq0$ and its absolute value increases
with the increase in the surface tension. It could also be observed
from the plots of Fig. 2 that, for this special case of simple shear
far-field load, the variation of the hoop stress with the angle $\theta$
is small with the maximum (minimum) achieved at $\theta=0$ ( $\theta=\pi/2$),
respectively. 

\begin{figure}[ht]
\centering\includegraphics[width=0.58\textwidth]{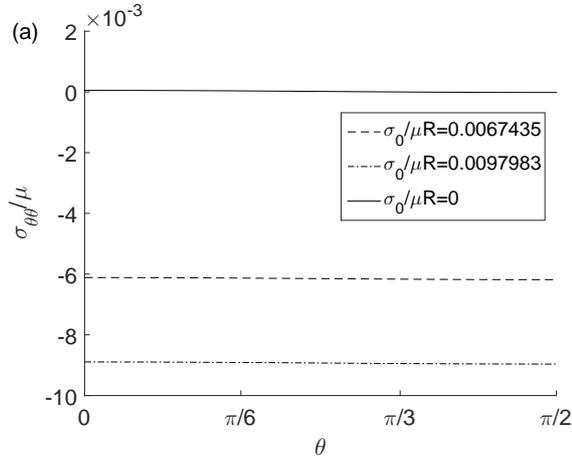} \caption{Normalized hoop stress $\sigma_{\theta\theta}/\mu$ along the cavity
surface, $\varphi=0$}
\label{fig_2} 
\end{figure}

The normalized stress $\sigma_{zz}/\mu$ along the same meridian line
$\varphi=0$ at the cavity surface is plotted on Fig. 3. Here too,
the surface tension has significant effect on the stress variation
reaching its maximum at $\theta=0$ and minimum at $\theta=\pi/2$.
We also notice that, for the simple shear far-field load, the variation
of $\sigma_{zz}/\mu$ with the angle $\theta$ is more pronounced
than that of $\sigma_{\theta\theta}/\mu$.

\begin{figure}[ht]
\centering\includegraphics[width=0.58\textwidth]{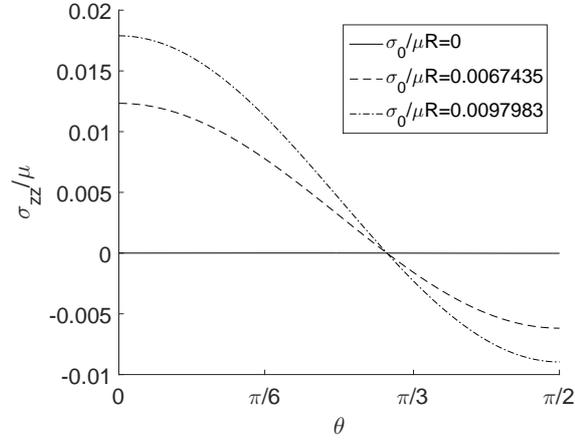} \caption{Normalized stress $\sigma_{zz}/\mu$ along the cavity surface, $\varphi=0$}
\label{fig_3} 
\end{figure}

\noindent To illustrate size-dependence of the surface related stresses,
consider an additional example of the matrix made from anodic alumina
($\mu=34.7$GPA and $\nu=0.3$) containing the cavity whose radius
varies from $R=5\,nm$ to $R=20\,nm.$ We fix the value of surface
tension as $\sigma_{0}=1.7N/m$ and that of the simple shear stress
as $\sigma_{d}=100$ MPA. The value of $R$-independent parameter
parameter $\gamma$ is chosen to be the same as the corresponding
expression in Eq. (\ref{eq:5.13}), while the remaining parameters
are chosen as

\begin{equation}
\mu_{0}=5.2321N/m,\:\lambda_{0}=10.4641N/m\label{eq:5.14}
\end{equation}
which for $R=5\,nm$ are consistent with the parameters used in expressions
of Eq. (\ref{eq:5.13}).

Fig. 4 illustrates the variation of the normalized hoop stress $\sigma_{\theta\theta}/\sigma_{d}$
along the line $\varphi=\pi/2$ at the cavity surface ($0\leq\theta\leq\pi/2$)
for various values of $R$. For comparison, we also plotted the variation
of the same but size-independent stress for the classical case (without
surface effects). It can be seen from the plots on Fig.4 that the
influence of surface effects on the normalized hoop stress diminishes
with the increase in $R$. 

\begin{figure}[ht]
\centering\includegraphics[width=0.58\textwidth]{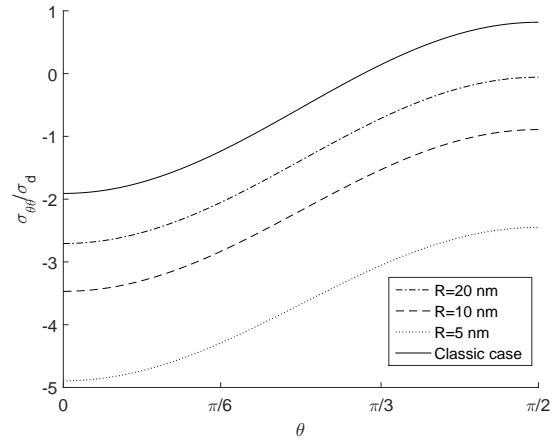} \caption{ Normalized hoop stress $\sigma_{\theta\theta}/\sigma_{d}$ along
the cavity surface, $\varphi=\pi/2$}
\label{fig_4} 
\end{figure}

\section{Effective properties of the isotropic particulate composites
with the complete Steigmann-Ogden model of interfaces}

\noindent The new analytical solution of previous section will now
be used to derive the single-inhomogeneity-based approximation formula
for the effective shear modulus of macroscopically isotropic composites
containing spherical inhomogeneities with the interfaces described
by the complete Steigmann-Ogden model. Similar formula for the effective
bulk modulus was reported in Zemlyanova and Mogilevskaya (2018a,b). 

To obtain shear modulus $\mu_{ef}$, we will use Maxwell's (1873)
concept of equivalent inhomogeneity, see e.g. McCartney (2010), Mogilevskaya
et al. (2012), and apply the same two-stage procedure as in Mogilevskaya
et al. (2018) and Zemlyanova and Mogilevskaya (2018a,b).

On the first stage, the problem involving a composite system containing
an inhomogeneity with the complete Steigmann-Ogden interface and subjected
to simple shear far-field load is solved twice, once with both the
external load and surface tension included and the second time with
just surface tension. The solution of the second problem is then subtracted
from that for the first problem. This procedure is needed to eliminate
the residual effects due to the presence of the surface tension, see
Mogilevskaya et al. (2010). As the results, the displacements in the
matrix will be given by Eqs. (\ref{eq:4.15})-(\ref{eq:4.18}) with
$D_{0}=0$ and remaining coefficients found from the system of Eqs.
(\ref{eq:C.1})-(\ref{eq:C.2}) of Appendix C. 

Then, we equate the coefficient $D_{4}$ for the leading $1/r^{2}$
(dipole) term in the obtained displacements with that related to the
problem involving the equivalent perfectly bonded inhomogeneity of
the same radius with the unknown shear modulus $\mu_{eq}$ to find
the value of that modulus.

On the second stage, the obtained modulus of the equivalent inhomogeneity
could be used in any single-inhomogeneity-based homogenization scheme
for perfectly bonded spherical particles. Here, we will use the following
Maxwell type approximation formulae, McCartney (2010):

\begin{align}
\frac{\mu_{ef}}{\mu} & =\frac{\mu{}_{eq}+\mu^{*}+c\mu^{*}\left(\mu_{eq}/\mu-1\right)}{\mu{}_{eq}+\mu^{*}-c\left(\mu_{eq}-\mu\right)}\label{eq:6.1}
\end{align}
with $c$ being the volume fraction of the equivalent inhomogeneity
and

\begin{equation}
\mu^{*}=\mu\frac{9K+8\mu}{6\left(K+2\mu\right)}=\frac{\mu}{2}\frac{9\lambda+14\mu}{3\lambda+8\mu}\label{eq:6.2}
\end{equation}

The solution for the problem of perfectly bonded inhomogeneity can
be extracted from Eq. (\ref{eq:4.18}) by assuming vanishing surface
parameters. The coefficient for the leading $1/r^{2}$ (dipole) term
for that case is $\left(5-4\nu\right)D_{4}^{eq}/\left(1-2\nu\right)$
in which

\begin{equation}
D_{4}^{eq}=-\frac{5}{2}\frac{\left(\mu_{eq}/\mu-1\right)\sigma_{d}^{\infty}R^{3}}{\left(9\lambda+14\mu\right)+2\left(\mu_{eq}/\mu\right)\left(8\lambda+3\mu\right)}\label{eq:6.7}
\end{equation}
with $\lambda=K^{\left(3\right)}-2\mu/3$.

Thus, assuming that $D_{4}^{eq}=D_{4}^{S-O}$, $\mu_{eq}$ could be
obtained as

\begin{equation}
\frac{\mu_{eq}}{\mu}=1-2R^{-3}\frac{4\mu+5\lambda}{1+4\left(D_{4}R^{-3}/\sigma_{d}^{\infty}\right)\left(3\mu+8\lambda\right)/5}\frac{D_{4}^{S-O}}{\sigma_{d}^{\infty}}\label{eq:6.8}
\end{equation}
in which $D_{4}^{S-O}/\sigma_{d}^{\infty}$ can be found from Eqs.
(\ref{eq:C.4}) of Appendix C.

Using the following notation:

\begin{equation}
\varLambda=\frac{2+5\lambda/2\mu}{\mu+4\mu\left(D_{4}^{S-O}R^{-3}/\sigma_{d}^{\infty}\right)\left(3\mu+8\lambda\right)/5}\frac{D_{4}^{S-O}R^{-3}}{\sigma_{d}^{\infty}}\label{eq:6.9}
\end{equation}
and substituting Eq. (\ref{eq:6.8}) into Eq. (\ref{eq:6.1}), we
obtain the following expression for the effective shear modulus $\mu_{ef}$:

\begin{equation}
\frac{\mu_{ef}}{\mu}=1-\frac{15c\left(\lambda+2\mu\right)\varLambda}{2\left(3\lambda+8\mu\right)\left[1-\left(1-c\right)\varLambda\right]+\left(9\lambda+14\mu\right)}\label{eq:6.10}
\end{equation}

To illustrate the effects of the Steigmann-Ogden interface parameters
on the effective shear modulus of the composite material with overall
isotropy, we consider the example of the preceding section with the
parameters given by Eqs. (\ref{eq:5.11})-(\ref{eq:5.13}) and use
Eq. (\ref{eq:6.10}) to estimate $\mu_{ef}$. Fig. 4 presents the
plots of the normalized effective shear modulus as a function of the
volume fraction of the cavities; the plot for the classical case with
$\mu_{0}=\lambda_{0}=\gamma=\sigma_{0}=0$ (solid line) is presented
there as well. It could be seen from Fig. 4 that while the effects
of surface elastic parameters lead to the increase in overall stiffness
of the composite, the influence of surface tension is insignificant
as its variation does not noticeable affect the value of $\mu_{ef}/\mu$.

\begin{figure}[ht]
\centering\includegraphics[clip,width=0.6\textwidth]{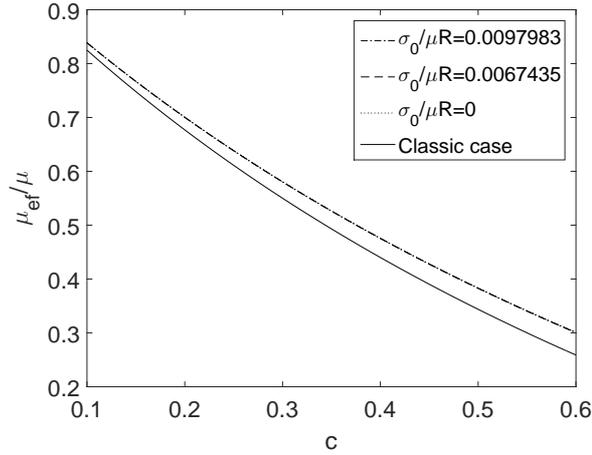} \caption{Normalized effective shear modulus for the composite material with
overall isotropy}
\label{fig_5} 
\end{figure}

\noindent To separate the influences of the \textit{elastic} interface
parameters, we assumed that $\sigma_{0}=0$ and evaluated, using the
single-inhomogeneity based estimate of Eq. (\ref{eq:6.10}), the values
of the normalized effective shear modulus for three different cases:
i) classical case with no surface effects, ii) the case with $\gamma/\mu=0$
(G-M model), and iii) the case with with $\gamma/\mu=0.00028382$
(S-O model). The remaining elastic parameters were taken from the
corresponding expressions of Eq. (\ref{eq:5.13}). The estimates for
all three cases are presented in the Table 1 for three values of the
volume fractions.

\begin{table}[h]
\centering \caption{Normalized shear modulus for the three interface models}
\vspace{0.3cm}
\begin{tabular}{|l|l|c|l|}
\hline 
\multicolumn{4}{|c|}{$\mu_{ef}/\mu$}\tabularnewline
\hline 
$c$  & Classic case  & G-M model  & S-O model \tabularnewline
\hline 
$0.1$  & $0.825$  & $0.838878$  & $0.838930$\tabularnewline
\hline 
$0.3$  & $0.550$  & $0.580938$  & $0.581057$\tabularnewline
\hline 
$0.5$  & $0.34375$  & $0.383570$  & $0.383725$\tabularnewline
\hline 
\end{tabular}
\end{table}

\noindent It could be seen from the table that, for the parameters
considered, the combined influence of the surface elasticity parameters
$\mu_{0},\:\lambda_{0}$ on the normalized effective shear modulus
is significantly larger than the influence of the bending parameter
$\gamma,$ especially for the smaller volume fractions of the cavities.

\section{Conclusions}

This paper presents two new important contributions to the studies
of the surface/interface effects in heterogeneous media. 

The first contribution is the analytical representations similar to
those of Christensen and Lo (1979) that could be used for solutions
of various problems involving spherical and circular material surfaces
that possess constant surface tension. These representations are derived
here for the first time. While in the present paper we used them for
the analysis of two specific interface models, we envision that they
could be used for the analysis of various Eshelby-type problems, including
new problems that might emerge in studies of surface effects. The
obtained representations significantly simplify the solutions of such
problems and can be useful in the comparative studies of various surface/interface
models and their influences on the local behavior of the fiber- and
particle reinforced nano composites. They should be especially valuable
for the researchers of engineering community as they are presented
in ready-to use forms and require relatively little algebra. In addition,
they could be extremely valuable to the researchers in the area of
Micromechanics as they allow for easy construction of various single-inhomogeneity-based
estimates for overall properties of such materials.

The second contribution is the new analytical solution of the problem
of a spherical inhomogeneity with the interface described by the \textit{complete}
Steigmann-Ogden model. This solution is important for the analysis
of the local and overall behavior of nano-scale materials and, at
the very least, the investigators who want to numerically solve more
complicated problems involving material surfaces could utilize our
solution as benchmark example. This solution was readily obtained
with the use of the new Christensen and Lo's type representations
for the displacements, while the alternative and may be more straightforward
approach to obtain it would require the use of tedious algebra of
spherical harmonics. The attractive feature of this new solution is
that all the coefficients involved in it are provided as closed-form
analytical expressions. This allows for accurate identification of
all governing dimensionless problem parameters and their influences
on the solution. 

\noindent \smallskip{}

\section*{Acknowledgements}

The first author (S.M.) gratefully acknowledges the support provided
by the Theodore W. Bennett Chair, University of Minnesota. The second
author (V.K.) acknowledges the support from Science and Technology
Center in Ukraine (STCU), award \#6247. The third author (A.Z.) gratefully
acknowledges the support from Simons Collaboration grant, award \#
319217. \smallskip{}

\appendix 
\section{The representations for the displacements
for the problem of Fig.1a}

The displacements in polar coordinates (Fig.1a) everywhere in the
composite system can be expressed in terms of two Kolosov-Muskhelishvili
potentials $\varphi(z)$ and $\psi(z)$ as, see Muskhelishvili (1959) 

\begin{equation}
2\mu\left[u_{r}(z)+iu_{\vartheta}(z)\right]=\exp\left(-i\vartheta\right)\left[\kappa\varphi(z)-z\overline{\varphi^{\prime}(z)}-\overline{\psi(z)}\right]\label{eq:3.3}
\end{equation}
in which the elastic parameters $\mu$ and $\kappa=3-4\nu$ are the
parameters of the specific phase: ($\mu,\,\kappa$) for the matrix
or ($\mu_{I},\,\kappa_{I}$) for the inhomogeneity.

In Mogilevskaya et al. (2008) and Zemlyanova and Mogilevskaya (2018),
it was shown that the potentials inside and outside of circular inhomogeneity
are

\textit{inside the inhomogeneity }($z=r\exp\left(i\vartheta\right),\:r<R$)

\begin{align}
\varphi(z) & =\frac{2\mu_{I}}{\kappa_{I}-1}ReA_{1}\frac{r}{R}\exp\left(i\vartheta\right)+\frac{2\mu_{I}}{\kappa_{I}}A_{3}\frac{r^{3}}{R^{3}}\exp\left(3i\vartheta\right)\label{eq:A.1}\\
\psi(z) & =-2\mu_{I}\left(\frac{3}{\kappa_{I}}A_{3}+\overline{A}_{-1}\right)\frac{r}{R}\exp\left(i\vartheta\right)\nonumber 
\end{align}

\textit{inside the} \textit{matrix }($z=r\exp\left(i\varphi\right),\:r>R$)

\begin{equation}
\varphi(z)=\frac{2}{\kappa+1}\left(-\omega^{\left(1\right)}A_{-1}+3\eta^{\left(2\right)}\overline{A}_{3}\right)\frac{R}{r}\exp\left(-i\vartheta\right)+\varphi^{\infty}(z)\label{eq:A.2}
\end{equation}

\begin{eqnarray*}
\psi(z) & = & \frac{2}{\kappa+1}\left\{ \left(\kappa-1\right)\left(\omega^{\left(0\right)}ReA_{1}+\frac{\sigma_{0}}{2}\right)\frac{R}{r}\exp\left(-i\vartheta\right)\right.\\
 & + & \left.\left[\left(\omega^{\left(1\right)}-\kappa\eta^{\left(2\right)}\right)A_{-1}+\omega^{\left(2\right)}\overline{A}_{3}\right]\frac{R^{3}}{r^{3}}\exp\left(-3i\vartheta\right)\right\} \\
 & + & \psi^{\infty}(z)
\end{eqnarray*}
in which 

\begin{eqnarray}
ReA_{1} & = & -\frac{1}{4\Delta_{1}}\sigma_{0}+\frac{\kappa+1}{4}\frac{1}{4\Delta_{1}}R\left(\sigma_{11}^{\infty}+\sigma_{22}^{\infty}\right)\label{eq:A.6}\\
A_{-1} & = & -\frac{\kappa+1}{4}\frac{\mu\kappa_{I}+\mu_{I}+3\kappa_{I}\eta^{\left(1\right)}}{\Delta_{2}}R\left(\sigma_{22}^{\infty}-\sigma_{11}^{\infty}-2i\sigma_{12}^{\infty}\right)\nonumber \\
A_{3} & = & -\frac{\kappa+1}{4}\frac{\kappa_{I}\eta^{\left(2\right)}}{\Delta_{2}}R\left(\sigma_{22}^{\infty}-\sigma_{11}^{\infty}+2i\sigma_{12}^{\infty}\right)\nonumber 
\end{eqnarray}

\begin{equation}
\varphi^{\infty}(z)=\frac{\sigma_{11}^{\infty}+\sigma_{22}^{\infty}}{4}z,\:\psi^{\infty}(z)=\frac{\sigma_{22}^{\infty}-\sigma_{11}^{\infty}+2i\sigma_{12}^{\infty}}{2}z\label{eq:A.7}
\end{equation}

and

\begin{eqnarray}
\Delta_{1} & = & \frac{\mu_{I}}{\kappa_{I}-1}+\frac{\mu}{2}+\eta\label{eq:A.8}\\
\Delta_{2} & = & \left(\mu+\kappa\mu_{I}\right)\left(\mu\kappa_{I}+\mu_{I}\right)+\eta^{\left(1\right)}\left[3\kappa_{I}\left(\mu+\kappa\mu_{I}\right)+\kappa\left(\mu\kappa_{I}+\mu_{I}\right)\right]\nonumber \\
 & + & 12\kappa\kappa_{I}\eta\left(\frac{\sigma_{0}}{4R}+\gamma\right)\nonumber 
\end{eqnarray}

Substituting the expressions of Eqs. (\ref{eq:A.1}), (\ref{eq:A.2})
into the representation of Eq. (\ref{eq:3.3}) (with the corresponding
elastic parameters) and performing some algebra, we obtain the representations
for the displacements in the polar coordinate system as

\textit{inside the inhomogeneity}

\begin{align}
u_{r}(z)+iu_{\vartheta}(z) & =ReA_{1}\frac{r}{R}+A_{3}\frac{r^{3}}{R^{3}}\exp\left(2i\vartheta\right)\label{eq:3.10}\\
+ & \left[\frac{3}{\kappa_{I}}\overline{A}_{3}\left(\frac{r}{R}-\frac{r^{3}}{R^{3}}\right)+A_{-1}\frac{r}{R}\right]\exp\left(-2i\vartheta\right)\nonumber 
\end{align}

\textit{inside the matrix}

\begin{align}
u_{r}(z)+iu_{\vartheta}(z) & =\frac{1}{\mu\left(\kappa+1\right)}\left\{ -\left(\kappa-1\right)\left(\omega^{\left(0\right)}ReA_{1}+\frac{\sigma_{0}}{2}\right)\frac{R}{r}\right.\nonumber \\
+ & \kappa\left(-\omega^{\left(1\right)}A_{-1}+3\eta^{\left(2\right)}\overline{A}_{3}\right)\frac{R}{r}\exp\left(-2i\vartheta\right)\nonumber \\
+ & \left(-\omega^{\left(1\right)}\overline{A}_{-1}+3\eta^{\left(2\right)}A_{3}\right)\frac{R}{r}\exp\left(2i\vartheta\right)\nonumber \\
+ & \left.\left[\left(-\omega^{\left(1\right)}+\kappa\eta^{\left(2\right)}\right)\overline{A}_{-1}-\omega^{\left(2\right)}A_{3}\right]\frac{R^{3}}{r^{3}}\exp\left(2i\vartheta\right)\right\} \nonumber \\
+ & u_{r}^{\infty}(z)+iu_{\vartheta}^{\infty}(z)\label{eq:3.11}
\end{align}
where $\exp\left(i\vartheta\right)=\cos\vartheta+i\sin\vartheta$

\begin{equation}
u_{r}^{\infty}(z)+iu_{\vartheta}^{\infty}(z)=\left(\kappa-1\right)\frac{\sigma_{11}^{\infty}+\sigma_{22}^{\infty}}{8\mu}r-\frac{\sigma_{22}^{\infty}-\sigma_{11}^{\infty}-2i\sigma_{12}^{\infty}}{4\mu}r\exp\left(-2i\vartheta\right)\label{eq:3.12}
\end{equation}

\begin{equation}
\omega^{\left(0\right)}=\frac{2\mu_{I}}{\kappa_{I}-1}-\frac{2\mu}{\kappa-1}+2\eta=K_{I}^{\left(2\right)}-K^{\left(2\right)}+2\eta\label{eq:3.13}
\end{equation}

\begin{equation}
\omega^{\left(1\right)}=\mu_{I}-\mu+\eta^{\left(1\right)}\label{eq:3.14}
\end{equation}

\begin{equation}
\omega^{\left(2\right)}=\mu_{I}\frac{\kappa}{\kappa_{I}}-\mu+3\kappa\eta^{\left(1\right)}+3\eta^{\left(2\right)}\label{eq:3.15}
\end{equation}

\begin{equation}
\eta^{\left(1\right)}=\eta+\gamma+\sigma_{0}/4R\label{eq:3.16}
\end{equation}

\begin{equation}
\eta^{\left(2\right)}=\eta-\gamma-\sigma_{0}/4R\label{eq:3.17}
\end{equation}

\begin{equation}
\eta=\frac{2\mu_{0}+\lambda_{0}}{4R}\label{eq:3.18}
\end{equation}
and the remaining parameters are defined by Eqs. (\ref{eq:A.6})-(\ref{eq:A.8})
.

From Eq. (\ref{eq:A.6}), it follows that, for the case of simple
shear far-field load ($\sigma_{11}^{\infty}=-\sigma_{22}^{\infty}=\sigma_{d}^{\infty},\:\sigma_{12}^{\infty}=0$),
the coefficients involved in Eqs. (\ref{eq:3.10}), (\ref{eq:3.11})
are

\begin{align}
ReA_{1} & =-\frac{\sigma_{0}}{2\left(K_{I}^{\left(2\right)}+\mu+2\eta\right)}\label{eq:3.20}\\
A_{-1} & =\frac{\kappa+1}{2}\frac{\mu\kappa_{I}+\mu_{I}+3\kappa_{I}\eta^{\left(1\right)}}{\Delta_{2}}R\sigma_{d}^{\infty}\nonumber \\
A_{3} & =\frac{\kappa+1}{2}\frac{\kappa_{I}\eta^{\left(2\right)}}{\Delta_{2}}R\sigma_{d}^{\infty}\nonumber 
\end{align}
where $\Delta_{2}$ is defined in Eq. (\ref{eq:A.8}).

\smallskip{}

\section{Vector representations for
the displacements for the problem of Fig.1b}

The following vector partial solutions are involved in representations
(\ref{eq:4.10}), (\ref{eq.:4.11}):

\begin{equation}
\begin{array}{ccc}
\mathbf{u}_{00}^{(3)}(\mathbf{x})=-2\left(1-2\nu\right)r\mathbf{S}_{00}^{(3)}/3, & \mathbf{u}_{22}^{(1)}(\mathbf{x})=\frac{r}{24}\left(\mathbf{S}_{22}^{(1)}+2\mathbf{S}_{22}^{(3)}\right), & \mathbf{u}_{22}^{(3)}(\mathbf{x})=\frac{r^{3}}{24}\left[\frac{7-4\nu}{21}\mathbf{S}_{22}^{(1)}+\frac{4\nu}{7}\mathbf{S}_{22}^{(3)}\right]\end{array}\label{eq:B.1}
\end{equation}

\begin{equation}
\begin{array}{ccc}
\mathbf{U}_{00}^{(1)}(\mathbf{x})=-\frac{1}{r^{2}}\mathbf{S}_{00}^{(3)}, & \mathbf{U}_{22}^{(1)}(\mathbf{x})=\frac{1}{r^{4}}\left(\mathbf{S}_{22}^{(1)}-3\mathbf{S}_{22}^{(3)}\right), & \mathbf{U}_{22}^{(3)}(\mathbf{x})=\frac{1}{3r^{2}}\left[\left(1-2\nu\right)\mathbf{S}_{22}^{(1)}+\left(5-4\nu\right)\mathbf{S}_{22}^{(3)}\right]\end{array}\label{eq:B.2}
\end{equation}
in which the vector spherical surface harmonics $\mathbf{S}_{ts}^{(i)}$
(Morse and Feshbach, 1953)
\begin{equation}
\mathbf{S}_{ts}^{(1)}=\mathbf{e}_{\theta}\frac{\partial}{\partial\theta}\chi_{t}^{s}+\frac{\mathbf{e}_{\varphi}}{\sin\theta}\frac{\partial}{\partial\varphi}\chi_{t}^{s},\quad\mathbf{S}_{ts}^{(2)}=\frac{\mathbf{e}_{\theta}}{\sin\theta}\frac{\partial}{\partial\varphi}\chi_{t}^{s}-\mathbf{e}_{\varphi}\frac{\partial}{\partial\theta}\chi_{t}^{s},\quad\mathbf{S}_{ts}^{(3)}=\mathbf{e}_{r}\chi_{t}^{s}.\label{eq:B.3}
\end{equation}
are expressed via the scalar surface harmonics $\chi_{t}^{s}(\theta,\varphi)$
as

\begin{equation}
\chi_{t}^{s}(\theta,\varphi)=P_{t}^{s}(\cos\theta)\exp\left(\mathrm{i}s\varphi\right)\label{eq:B.4}
\end{equation}
and $P_{t}^{s}(\cos\theta)$ are associated Legendre polynomials.

The constants $D_{0}$, $D_{3}$ $D_{4}$, $A_{0}$, $A_{1}$, $A_{2}$
of Eqs. (\ref{eq:4.12})-(\ref{eq:4.17}) are found from the system
of linear equations that represent the interface and far-field conditions
(Kushch et al., 2011, Kushch, 2013). In our notations, they are
\begin{itemize}
\item continuity of displacements 
\begin{gather}
\frac{1}{R^{3}}D_{0}=A_{0}\label{eq:B.7}\\
\frac{\sigma_{d}^{\infty}}{4\mu}-\frac{2}{R^{5}}D_{3}+\frac{2}{R^{3}}D_{4}=A_{1}-R^{2}\frac{7-4\nu_{I}}{1-2\nu_{I}}A_{2}\nonumber \\
\frac{\sigma_{d}^{\infty}}{4\mu}+\frac{3}{R^{5}}D_{3}+\frac{1}{R^{3}}\frac{5-4\nu}{1-2\nu}D_{4}=A_{1}-R^{2}\frac{6\nu_{I}}{1-2\nu_{I}}A_{2}\nonumber 
\end{gather}
\item interface conditions of Eqs. (\ref{eq:4.1})-(\ref{eq:4.3}) 
\begin{gather}
\frac{2}{R^{3}}D_{0}=-\frac{\mu_{I}}{\mu}\frac{\left(1+\nu_{I}\right)}{(1-2\nu_{I})}A_{0}-f_{0}\label{eq:B.8}\\
\frac{\sigma_{d}^{\infty}}{4\mu}+\frac{2}{R^{3}}\left(\frac{4}{R^{2}}D_{3}+\frac{1+\nu}{1-2\nu}D_{4}\right)=\frac{\mu_{I}}{\mu}\left(A_{1}-R^{2}\frac{7+2\nu_{I}}{1-2\nu_{I}}A_{2}\right)-6f_{1}\nonumber \\
\frac{\sigma_{d}^{\infty}}{4\mu}+\frac{2}{R^{3}}\left(-\frac{6}{R^{2}}D_{3}+\frac{\nu-5}{1-2\nu}D_{4}\right)=\frac{\mu_{I}}{\mu}\left(A_{1}+R^{2}\frac{3\nu_{I}}{1-2\nu_{I}}A_{2}\right)-3f_{2}\nonumber 
\end{gather}
In Eq. (\ref{eq:B.8}), 
\begin{gather}
f_{0}=-\frac{\sigma_{0}}{\mu R}+\frac{\eta_{0}}{\mu R}A_{0}\label{eq:B.9}\\
f_{1}=\frac{(\mu_{0}-\sigma_{0})-3(\lambda_{0}+2\mu_{0})}{6\mu R}\left(A_{1}-R^{2}\frac{7-4\nu_{I}}{1-2\nu_{I}}A_{2}\right)\nonumber \\
+\frac{\lambda_{0}+\mu_{0}+\sigma_{0}}{3\mu R}\left(A_{1}-R^{2}\frac{6\nu_{I}}{1-2\nu_{I}}A_{2}\right)\nonumber \\
f_{2}=\frac{\lambda_{0}+\mu_{0}+\sigma_{0}}{\mu R}\left(A_{1}-R^{2}\frac{7-4\nu_{I}}{1-2\nu_{I}}A_{2}\right)\nonumber \\
-2\frac{\lambda_{0}+\mu_{0}+2\sigma_{0}}{3\mu R}\left(A_{1}-R^{2}\frac{6\nu_{I}}{1-2\nu_{I}}A_{2}\right)\nonumber 
\end{gather}
\end{itemize}
in which $\eta_{0}$ is defined by Eq. (\ref{eq:4.23}).

After some algebraic manipulations with Eqs. (\ref{eq:B.7})-(\ref{eq:B.9}),
one can arrive to the following expressions for the coefficients $A_{0},\:D_{0}$:
\begin{equation}
A_{0}=D_{0}/R^{3}=\frac{2\sigma_{0}/R}{4\mu+3K_{I}^{\left(3\right)}+2\eta_{0}}\label{eq:B.10}
\end{equation}

The analytical expressions for the remaining coefficients will be
presented in Appendix C as a particular case of more general expressions
for the coefficients related to the Steigmann-Ogden model.

\section{Coefficients involved in the
representations of Eqs. (\ref{eq:2.4}), (\ref{eq:2.5}) }

The continuity conditions for the displacements produce two linear
equations for the unknown coefficients $A_{1}$, $A_{2}$, $D_{3}$,
$D_{4}$ in the representations of Eqs. (\ref{eq:2.4}), (\ref{eq:2.5}),
(\ref{eq:4.12})-(\ref{eq:4.17}): 
\begin{equation}
A_{1}R-3\frac{\lambda_{I}}{\mu_{I}}A_{2}R^{3}-3D_{3}R^{-4}-\frac{3\lambda+5\mu}{\mu}D_{4}R^{-2}=D_{1}R\label{eq:C.1}
\end{equation}
\[
A_{1}R-\frac{5\lambda_{I}+7\mu_{I}}{\mu_{I}}A_{2}R^{3}+2D_{3}R^{-4}-2D_{4}R^{-2}=D_{1}R
\]

Substitution of the representations of Eqs. (\ref{eq:2.4}), (\ref{eq:2.5}),
(\ref{eq:4.12})-(\ref{eq:4.17}) into Eqs. (\ref{eq:5.1})-(\ref{eq:5.3})
results in the following two additional linear equations for the unknown
constants $A_{1}$, $A_{2}$, $D_{3}$, $D_{4}$: 
\begin{equation}
C_{31}A_{1}R+C_{32}A_{2}R^{3}-24\mu D_{3}R^{-4}-(18\lambda+20\mu)D_{4}R^{-2}=-2\mu D_{1}R\label{eq:C.2}
\end{equation}
\[
C_{41}A_{1}R+C_{42}A_{2}R^{3}+8\mu D_{3}R^{-4}+(3\lambda+2\mu)D_{4}R^{-2}=-\mu D_{1}R
\]
where the parameters $C_{31}$, $C_{32}$, $C_{41}$, $C_{42}$ are
\[
C_{31}=-2\mu_{I}+\frac{2}{R}(\lambda_{0}+\mu_{0}-\sigma_{0})-6\gamma
\]
\begin{equation}
C_{32}=-3\lambda_{1}-\frac{6}{R}(\lambda_{0}+\mu_{0})\frac{3\lambda_{I}+7\mu_{I}}{\mu_{I}}-\frac{6}{R}\sigma_{0}\frac{\lambda_{I}+7\mu_{I}}{\mu_{I}}+6\gamma\frac{\lambda_{I}-7\mu_{I}}{\mu_{I}}\label{eq:C.3}
\end{equation}
\[
C_{41}=-\mu_{I}-\frac{1}{R}(3\mu_{0}+\lambda_{0}-\sigma_{0})+\gamma
\]
\[
C_{42}=8\lambda_{I}+7\mu_{1}+\frac{\mu_{0}}{R}\frac{19\lambda_{I}+35\mu_{I}}{\mu_{I}}+3\frac{\lambda_{0}}{R}\frac{3\lambda_{I}+7\mu_{I}}{\mu_{I}}-\frac{\sigma_{0}}{R}\frac{\lambda_{I}-7\mu_{I}}{\mu_{I}}-\gamma\frac{\lambda_{I}-7\mu_{I}}{\mu_{I}}
\]

As the Steigmann-Ogden model reduces to that of Gurtin-Murdoch, Eqs.
(\ref{eq:C.1})-(\ref{eq:C.3}) should reproduce the corresponding
equations of Eqs. (\ref{eq:B.7})-(\ref{eq:B.9}), when the bending
parameters vanish, i.e. $\gamma=0.$ After some algebraic manipulations,
not shown here, we have verified that this was, indeed, the case.

The system of Eqs. (\ref{eq:C.1})-(\ref{eq:C.2}) can be solved analytically
leading to the following expressions for the coefficients $A_{1}$,
$A_{2}$, $D_{3}$, $D_{4}$ : 
\[
A_{1}/D_{1}=F(E_{22}-E_{12})/(E_{11}E_{22}-E_{12}E_{21})
\]
\begin{equation}
A_{2}/D_{1}=F(E_{11}-E_{21})R^{-2}/(E_{11}E_{22}-E_{12}E_{21})\label{eq:C.4}
\end{equation}
\[
D_{3}/D_{1}=-\left[\left(A_{1}/D_{1}\right)\left((3\lambda+2\mu)C_{31}+(18\lambda+20\mu)C_{41}\right)R^{5}\right.
\]
\[
\left.+\left(A_{2}/D_{1}\right)\left((3\lambda+2\mu)C_{32}+(18\lambda+20\mu)C_{42}\right)R^{8}+24\mu(\mu+\lambda)R^{5}\right](8\mu)^{-1}(9\lambda+14\mu)^{-1}
\]
\[
D_{4}/D_{1}=\left[\left(A_{1}/D_{1}\right)(C_{31}+3C_{41})R^{3}+\left(A_{2}/D_{1}\right)(C_{32}+3C_{42})R^{5}+5\mu R^{3}\right](9\lambda+14\mu)^{-1}
\]
where the parameters $E_{11}$, $E_{12}$, $E_{21}$, $E_{22}$, $F$
are given by 
\[
E_{11}=1-\left[(3\lambda+10\mu)C_{31}+(18\lambda+44\mu)C_{41}\right](4\mu)^{-1}(9\lambda+14\mu)^{-1}
\]
\begin{equation}
E_{12}=-\left[(3\lambda+10\mu)C_{32}+(18\lambda+44\mu)C_{42}\right](4\mu)^{-1}(9\lambda+14\mu)^{-1}-(5\lambda_{1}+7\mu_{1})/\mu_{1}\label{eq:C.5}
\end{equation}
\[
E_{21}=1-\left[(15\lambda+34\mu)C_{31}+6(3\lambda+10\mu)C_{41}\right](8\mu)^{-1}(9\lambda+14\mu)^{-1}
\]
\[
E_{22}=-3\lambda_{1}/\mu_{1}-\left[(15\lambda+34\mu)C_{32}+6(3\lambda+10\mu)C_{42}\right](8\mu)^{-1}(9\lambda+14\mu)^{-1}
\]
\[
F=15(\lambda+2\mu)/(9\lambda+14\mu)
\]

\begin{center}
\textbf{List of Figures}
\par\end{center}

Fig. 1. a) Circular and b) spherical inhomogeneity in an infinite
matrix 

Fig. 2. Normalized hoop stress $\sigma_{\theta\theta}/\mu$ along
the cavity surface, $\varphi=0$

Fig. 3. Normalized stress $\sigma_{zz}/\mu$ along the cavity surface,
$\varphi=0$

Fig. 4. Normalized hoop stress $\sigma_{\theta\theta}/\sigma_{d}$
along the cavity surface, $\varphi=\pi/2$

Fig. 5. Normalized effective shear modulus for the composite material
with overall isotropy\medskip{}

\begin{center}
\textbf{List of Tables}
\par\end{center}

Table 1. Normalized shear modulus for the three interface models

\vspace{0.2cm}

\textbf{References}

\begin{thebibliography}{100}
\bibitem{Benveniste1989}Y. Benveniste, G.J. Dvorak, T. Chen. Stress
fields in composites with coated inclusions. Mech. Mater. 7, (1989), 305-317.

\bibitem{Cahn1982}J.W. Cahn, F. Larché. Surface stress and the
chemical equilibrium of small crystals\textemdash II. Solid particles
embedded in a solid matrix. Acta Metall. 30, (1982), 51-56.

\bibitem{Chatzigeorgiou2017} G. Chatzigeorgiou, F. Meraghni,  A. Javili.
Generalized interfacial energy and size effects in composites. J.
Mech. Phys. Solids 106, (2017), 257-282.

\bibitem{Chhapadia2011} P. Chhapadia,  P. Mohammadi,  P. Sharma. Curvature-dependent
surface energy and implications for nanostructures. J. Mech. Phys.
Solids 59, (2011), 2103-2115.

\bibitem{Christensen1979}R.M. Christensen, K.H. Lo. Solutions for effective
shear properties in three phase sphere and cylinder models. J. Mech.
Phys. Solids 27, (1979), 315-330.

\bibitem{Dingreville2005} R. Dingreville, J.M. Qu, M. Cherkaoui. Surface
free energy and its effect on the elastic behavior of nano-sized particles.
J. Mech. Phys. Solids 53, (2005), 1827-1854.

\bibitem{Duan2005}H.L. Duan,  J. Wang, Z.P. Huang, B.L. Karihaloo. Size-dependent effective elastic constants of solids containing
nano-inhomogeneities with interface stress. J. Mech. Phys. Solids
53, (2005), 1574-1596.

\bibitem{Duan2006}H.L. Duan, J. Wang, B.L. Karihaloo, Z.P. Huang. Nanoporous materials can be made stiffer that non-porous counterparts
by surface modification. Acta Materialia 54, (2006), 2983-2990.

\bibitem{Duan2007} H.L.Duan, X. Yi, Z.P. Huang, J. Wang.
A united scheme for prediction of effective moduli of multiphase composites
with interface effects. Part I: Theoretical framework. Mech. Mater.
39, (2007), 81-93.

\bibitem{Goodier1933} J.N. Goodier. Concentration of stress around
spherical and cylindrical inclusions and flaws. J. Appl. Mech. 55, (1933), 39-44.

\bibitem{Gurtin1975}M.E. Gurtin, A.I. Murdoch. A continuum theory
of elastic material surfaces. Arch. Ration. Mech. Anal. 57, (1975), 291-323.

\bibitem{gurtin-Murd2} M.E. Gurtin, A.I. Murdoch. Surface
stress in solids. Int. J. Solid. Struct. 14, (1978), 431-440.

\bibitem{Han2018} Z. Han, S.G. Mogilevskaya, D. Schillinger.
Local fields and overall transverse properties of unidirectional composite
materials with multiple nanofibers and Steigmann-Ogden interfaces.
Int. J. Solid. Struct. 147, (2018), 166-182.

\bibitem{Hashin1991} Z. Hashin. The spherical inclusion with imperfect
interface. J. Appl. Mech. 58, (1991), 444-449.

\bibitem{He2006}L.H. He, Z.R. Li. Impact of surface stress
on stress concentration. Int. J. Solid. Struct. 43, (2006), 6208-6219.

\bibitem{Herve1993} E. Herve,  A. Zaoui. n-Layered inclusion-based
micromechanical modelling. Int. J. Engng. Sci. 31, (1993), l-10.

\bibitem{Huang1993}Huang, Z.P., Wang, J., 2006. A theory of hyperelasticity
of multi-phase media with surface/interface energy effect. Acta Mech.
182, (1993), 195-210.

\bibitem{Jammes2009} M. Jammes, S.G. Mogilevskaya, S.L. Crouch.
Multiple circular nano-inhomogeneities and/or nano-pores in one of
two joined isotropic elastic half-planes. Eng. Anal. Bound. Elem.
33, (2009), 233-248.

\bibitem{Javili2017} A. Javili, N.S. Ottosen, M. Ristinmaa, J. Mosler. Aspects of interface elasticity theory. Math. Mech. Solids
23, (2017), 1004-1024.

\bibitem{Kushch2013} V.I. Kushch. Micromechanics of composites:
Multipole expansion approach. Elsevier, (2013).

\bibitem{Kushch2018} V.I. Kushch. Stress field and effective elastic
moduli of periodic spheroidal particle composite with Gurtin-Murdoch
interface. Int. J. Eng. Sci. 132, (2018), 79-96.

\bibitem{Kushch2014}  V.I. Kushch, V.S. Chernobai, G.S. Mishuris.
Longitudinal shear of a composite with elliptic nanofibers: local
stresses and effective stiffness. Int. J. Eng. Sci. 84, (2014), 79-94. 

\bibitem{Kushch2011} V.I. Kushch, S.G. Mogilevskaya, H.K. Stolarski,  S.L. Crouch. Elastic interaction of spherical nanoinhomogeneities
with Gurtin\textendash Murdoch type interfaces. J. Mech. Phys. Solids
59, (2011), 1702=1716.

\bibitem{Kushch2013a}  V.I. Kushch, S.G. Mogilevskaya, H.K. Stolarski. Elastic fields and effective moduli of particulate nanocomposites
with the Gurtin-Murdoch model of interfaces. Int. J. Solid. Struct.
50, (2013), 1141-1153.

\bibitem{Kushch2018a} V.I. Kushch, S.V. Shmegera, V.V. Mykhas'kiv. Multiple spheroidal cavities with surface stress as a model
of nanoporous solid. Int. J. Solid. Struct. 152\textendash 153, (2018), 261-271.

\bibitem{Lim2006} C.W. Lim, Z.R. Li, L.H. He. Size dependent,
non-uniform elastic field inside a nano-scale spherical inclusion
due to interface stress. Int. J. Solid. Struct. 43, (2006), 5055-5065.

\bibitem{Love1927} A.E.H. Love. A Treatise on the Mathematical
Theory of of Elasticity (4th edn). University Press, Cambridge, (1927).

\bibitem{Lurie2005}  A.I. Lurie. Theory of Elasticity. Springer,
New York (English translation of Russian 1970 edition, Nauka, Moscow), (2005).

\bibitem{Maxwell}J. C. Maxwell. Treatise on Electricity and
Magnetism, vol. 1, Clarendon Press, Oxford (3rd edn), (1892).

\bibitem{McCartney2010}L.N. McCartney. Maxwell's far-field
methodology predicting elastic properties of multi-phase composites
reinforced with aligned transversely isotropic spheroids. Philos.
Mag. 90, (2010), 4175-4207.

\bibitem{Mi2006}C. Mi, D.A. Kouris. Nanoparticles under the
influence of surface/interface elasticity. Mech. Mater. Struct. 1, (2006), 763-791.

\bibitem{Miller2000}R.E. Miller, V.B. Shenoy. Size-dependent elastic
properties of nanosized structural elements. Nanotechnology 11, (2000), 139-147.

\bibitem{Mogilevskaya2010} S.G. Mogilevskaya, S.L. Crouch, A. LaGrotta, H.K. Stolarski. The effects of surface elasticity and surface tension
on the transverse overall elastic behavior of unidirectional nano-composites.
Compos. Sci. 70, (2010), 427-434.

\bibitem{Mogilevskaya2008} S.G. Mogilevskaya, S.L. Crouch, H.K. Stolarski. Multiple interacting circular nano-inhomogeneities with
surface/interface effects. J. Mech. Phys. Solids 56, (2008), 2298-2327.

\bibitem{M-Maxwell}S.G. Mogilevskaya, H.K. Stolarski, S.L. Crouch. On Maxwell\textquoteright s concept of equivalent inhomogeneity:
when do the interactions matter? J. Mech. Phys. Solids 60, (2012), 391-417.

\bibitem{Mogilevskaya2018} S.G. Mogilevskaya, A.Y. Zemlyanova, M. Zammarchi. On the elastic far-field response of a two-dimensional coated
circular inhomogeneity: Analysis and applications. Int. J. Solid.
Struct. 130-131, (2018), 199-210.

\bibitem{Morse1953}P.M. Morse, H. Feshbach. In: Methods of Theoretical Physics. McGraw-Hill, NewYork, (1953).

\bibitem{Muskhelishvili1959}N.I. Muskhelishvili. Some Basic Problems of
the Mathematical Theory of Elasticity. Noordhoff, Groningen, 1959.

\bibitem{Ru2010}C.Q. Ru. Simple geometrical explanation of
Gurtin-Murdoch model of surface elasticity with clarification of its related versions. Sci. China Phys., Mech. \& Astron. 53, (2010), 536-544.

\bibitem{Savin1961}G.N. Savin. Stress concentrations around holes.
Pergamon Press, Oxford, (1961).

\bibitem{Sharma2002}P. Sharma, S. Ganti. Interfacial elasticity
corrections to size-dependent strain-state of embedded quantum dots. Phys. Status Solidi 234, (2002), R10-R12.

\bibitem{sharma-Ganti-2004} P. Sharma, S. Ganti. Size-dependent
Eshelby's tensor for embedded nano-inclusions incorporating surface/interface
energies. J. Appl. Mech. 71, (2004), 663-671.

\bibitem{sharma-et-al} P. Sharma, S. Ganti, N. Bhate. Effect
of surfaces on the size-dependent elastic state of nano-inhomogeneities.
Appl. Phys. Lett. 82, (2003), 535-537.

\bibitem{SteigmannOgden1997} D.J. Steigmann, R.W. Ogden.
Plain deformations of elastic solids with intrinsic boundary elasticity.
Proc. R. Soc. London A 453, (1997), 853-877.

\bibitem{Steigman} D.J. Steigmann, R.W. Ogden. Elastic surface-substrate
interactions. Proc. R. Soc. London A 455, (1999), 437-474.

\bibitem{Xu2016} Y. Xu, Q-C. He, S-T. Gu. Effective elastic
moduli of fiber-reinforced composites with interfacial displacement
and stress jumps. Int. J. Solids Struct. 80, (2016), 146-157.

\bibitem{Zemlyanova2018a} A.Y. Zemlyanova, S.G. Mogilevskaya. On spherical
inhomogeneity with Steigmann\textendash Ogden interface. J. Appl.
Mech. 85, (2018a), 121009-121009-10.

\bibitem{Zemlyanova2018b} A.Y. Zemlyanova, S.G. Mogilevskaya. Circular
inhomogeneity with Steigmann\textendash Ogden interface: Local fields,
neutrality, and Maxwell\textquoteright s type approximation formula.
Int. J. Solids Struct. 135, (2018b), 85-98.


\end{thebibliography}
\end{document}